\documentclass
{amsart}
\usepackage[utf8]{inputenc}
\usepackage{color}
\usepackage{graphicx}
\usepackage{floatrow}
\usepackage{comment}
\usepackage{varioref}
\usepackage{hyperref}
\usepackage{cleveref}
\usepackage{setspace}
\usepackage{lineno}

\newtheorem{proposition}{Proposition}
\newtheorem{conjecture}{Conjecture}

\newcommand{\puiseux}[2]{#1\{\!\{#2\}\!\}}

\author{Anders Jensen}
\address{{\tt jensen@math.au.dk}, Department of Mathematics, Aarhus University}
\author{Anton Leykin}
\address{{\tt leykin@math.gatech.edu}, School of Mathematics, Georgia Tech}

\newcommand{\ZZ}{{\mathbb Z}}
\newcommand{\QQ}{{\mathbb Q}}

\newcommand{\CC}{{\mathbb C}}
\newcommand{\VV}{{\mathbb V}}

\newcommand{\constC}{C}

\newcommand{\masses}{{m_1,\dots,m_n}}

\begin{document}
\title{Smale's 6th problem for generic masses}
\begin{abstract}
    We give a new method to attempt to prove that, for a given $n$, there are finitely many equivalence classes of planar central configurations in the Newtonian $n$-body problem for generic masses. 
    
    The human part of the proof relies on tropical geometry. The crux of our technique is in a computation that we have completed for $n\leq 5$, thus confirming the celebrated result of Albouy and Kaloshin. 
\end{abstract}
\maketitle

\section{Introduction}

All central configurations for $n=3$ were found by Euler and Lagrange explicitly. For $n=4$, Hampton and Moeckel~\cite{hampton2006finiteness} give an upper bound, 8472, on the number of solutions to a system of sparse polynomial equations describing the problem. Albouy and Kaloshin in~\cite{albouy2012finiteness} prove that the number of planar central configurations is finite for the five bodies of \emph{generic} masses. Their article gives a beautiful introduction to the history of the problem, which we shall not repeat, and an intricate proof of finiteness for $n=5$ in case of generic masses, an \emph{alternative} to which we develop.  

Our approach rests on tropical algebraic geometry apparatus.  
The approaches of~\cite{hampton2006finiteness} and \cite{albouy2012finiteness} can, in fact, be also described in the language of tropical geometry, which captures asymptotic behavior of algebraic solutions. 
The work of \cite{HJ11} combines polyhedral and polynomial computations to obtain equations for the set of possible exceptional choices of masses for which finiteness of central configurations in the spatial 5-body problem might not hold. In our new approach for showing generic finiteness in the planar 5-body problem we consider masses in a larger field with valuation, which empowers us with more tropical tools. This new technique leads us to the result using only polyhedral computations: in contrast with ~\cite{albouy2012finiteness,HJ11,hampton2006finiteness}, no other polynomials than those in the original equations are derived.

We believe that our approach could be fruitful in attacking the case of $n=6$ generic masses, albeit our computations have not been successful so far; see \Cref{sec:n=6}. Counterintuitively, more special open cases for $n=5$ involving either symmetry or additional relations on masses may be harder given the nature of our method; see \Cref{sec:special}. 

\section{Description of the problem}

\subsection{Newtonian mechanics}
We say that $n$ point masses $m_1,\dots,m_n$ with distinct positions $p_i=(x_i,y_i)$ 
form a \emph{central configuration} if the following equations are satisfied:
\begin{equation}\label{eq:main-xy}
  \constC p_i = \sum_{j\neq i} m_j \frac{p_j-p_i}{\|p_j-p_i\|^3} 
  \quad\quad (i\in 1,\dots,n) 
\end{equation}
for a nonzero constant $\constC$.

One can give the following physical meaning
to a configuration of the bodies satisfying \cref{eq:main-xy} (as we speak of physics all symbols take real values; masses and $\constC$ are positive).

Note that, if the $i$-th equation is multiplied by $m_i$, the right hand side equals (up to a constant multiplier) the Newton's \emph{gravitational force} acting on the $i$-th body. Now one can imagine the bodies rotating around the origin with a certain (same) constant angular velocity such that the left hand side equals the negative of Newton's \emph{centrifugal force} pulling the $i$-th body away from the origin.
To summarize, the bodies remain in the same relative positions to each other (as every body rotates around the origin in absolute coordinates).     

Observe that scaling a solution to \eqref{eq:main-xy} --- multiplying all coordinates by a constant --- gives a solution to the same system of equations for some value of constant $\constC$. 

Also observe a rotational symmetry: given a solution, rotating every body by the same angle produces another solution.

\subsection{Equations}\label{sec:AC}

The symmetric and nonsymmetric Albouy-Chenciner equations appearing in this section stem from~\cite{albouy1997probleme}. 
The derivation of the symmetric equations are carefully explained in~\cite{hampton2006finiteness}.
Roberts observed that symmetrizing is unnecessary in the explanation, thereby obtaining the nonsymmetric equations~\cite{Gareth}

One advantage of the following formulation over its alternatives --- for instance, the one in~\cite{albouy2012finiteness} --- is that it uses only distances, which are rotationally invariant. A certain normalization is incorporated making sure that, with respect to scaling, one representative per orbit is taken.

We consider the polynomial ring in $\binom{n}{2}$ variables, $r_{ij}$ with $1\leq i < j\leq n$, the pairwise distances among $n$ bodies and set the convention: $r_{ji} := r_{ij}$ for $j>i$ and $r_{ii}:=0$ for all $1\leq i \leq n$.

We use the following equations: 
\begin{itemize}
    \item ${(n-1)n}$ {\bf Albouy-Chenciner equations}:
    
    Let $S_{i,j}:=r_{ij}^{-3}-1$ for $i\not=j$ and $S_{i,i}:=0$. It was observed that the following polynomial vanishes at normalised distance vectors of central configurations:
    $$ g_{i,j}=\sum_k m_kS_{i,k}(r_{jk}^2-r_{ik}^2-r_{ij}^2).$$
    
    \item $\binom{n}{4}$ {\bf two-dimensional Cayley-Menger determinants}: 
    
    These ensure that the configuration is planar: for any choice of bodies $i_1,i_2,i_3,i_4$ the three-dimensional volume of their convex hull must be zero.
$$
\textup{det}\left(\begin{bmatrix}
0 & 1 & 1 & 1 & 1 \\
1 & 0 & r_{i_1i_2}^2 & r_{i_1i_3}^2 & r_{i_1i_4}^2 \\
1 & r_{i_1i_2}^2 & 0 & r_{i_2i_3}^2 & r_{i_2i_4}^2 \\
1 & r_{i_1i_3}^2 & r_{i_2i_3}^2 & 0 & r_{i_3i_4}^2 \\
1 & r_{i_1i_4}^2 & r_{i_2i_4}^2 & r_{i_3i_4}^2 & 0 
\end{bmatrix}\right)=0 
 $$
    \item $\frac{(n-1)n}{2}$ {\bf symmetric Albouy-Chenciner equations}:
    $$ f_{i,j}=g_{i,j}+g_{j,i}.$$
    
    
\end{itemize}

\subsection{Formulation of Smale's 6th problem and the conjecture for generic masses}

The original statement of Smale's problem --- one of the 18 ``problems for the next century''~\cite{smale1998mathematical} --- is stated with \emph{real} central configurations in mind:
\begin{quote} \em
    Given positive real numbers $m_l,\dots,m_n$ as the masses
in the $n$-body problem of celestial mechanics, is the
number of relative equilibria finite? 
\end{quote}

We aim to show that the assertion holds for \emph{almost all} masses.  
The masses are not assumed to be positive real numbers.
\begin{conjecture}\label{conj:generic}
For \emph{generic} values of masses, the number of normalized central configurations in the $n$-body problem over any field (of characteristic 0) is finite.
\end{conjecture}
One exact meaning of ``generic'' is that the conclusion holds for a complement of a hypersurface (described by a polynomial in $m_1,\dots,m_n$) in the parameter space (space of masses).

\subsection{Generic dimension}
Our main system is a system of polynomials with coefficients that are rational functions of $\masses$, therefore, one can attack \Cref{conj:generic} trying to confirm that the dimension of the solution set is $0$ generically over \emph{any} field $K\supseteq\QQ$. If successful, the result will hold for \emph{every} field (of characteristic $0$).\; Indeed, the locus of points in the space of masses for which the dimension is not $0$ would have to be contained in a proper subset cut out by polynomials with coefficients over $\QQ$, the original field of definition.

Warning: the tropical geometry terms in the paragraph below are defined in the (next) \Cref{sec:TG} and the full computational approach is described in \Cref{sec:computations}. 

Our strategy will be to use the field of Puiseux series as $K$ and consider the tropical variety of the given problem. One key fact is that the fiber of a point under the tropicalization map is Zariski dense. Carrying out the computations for \emph{one} tropical point (fixed valuations of masses) and showing that the dimension is $0$ will mean that the dimension is $0$ for almost \emph{every} point (values of masses in $K$) and, therefore, confirm \Cref{conj:generic}.    


\section{Basic principles of tropical geometry}\label{sec:TG}

We shall attack \Cref{conj:generic} over $K=\puiseux{\CC}{t}$, the field of \emph{Puiseux series}, a power series with rational exponents having a common denominator 
$$
\sum_{i=i_0}^{\infty} c_it^{i/D}, \quad i_0 \in \ZZ,\ D\in\ZZ_{>0},\ c_i\in\CC,\ c_{i_0}\neq 0.
$$ 
The \emph{valuation} map $K^* \rightarrow \QQ$, taking a series in $K^*:=K\setminus\{0\}$ to its lowest exponent, is applied coordinatewise to give the \emph{tropicalization} map $T: (K^*)^N \to \QQ^N$. 

The field $K$ is algebraically closed, which allows us to define a \emph{variety} as a common set of zeroes of a set of polynomials. For a system of polynomials $F$, the variety it cuts out is denoted by $\VV(F)$.

\subsection{Dimension is preserved.}
Tropical geometry may be seen as an ``asymptotic'' shadow of algebraic geometry. While tropicalization $T$ forgets some information about a variety $X$, some invariants are still captured by the variety's behavior ``at infinity'' and can be harvested from $T(X)$.    

For a variety $X \subset (K^*)^N$, define its \emph{tropical variety} as $T(X) \subset \QQ^N$. A tropical variety, in fact, may be viewed as a (finite) polyhedral complex. The dimension of $T(X)$ is the largest dimension of the polyhedra it contains. By Bieri-Groves Theorem \cite[Theorem~9.6]{groves1984geometry}, $\dim T(X) = \dim X$.

\subsection{Balancing condition}
A \emph{tropical curve}, a 1-dimensional tropical variety, is a finite complex of points, line segments, and halflines. It is \emph{balanced} locally, i.e., at every point the direction vectors (scaled appropriately\footnote{An appropriate scaling can be defined precisely, but is not necessary for our purposes.}) of incident line segments and halflines sum to zero.   A consequence of balancing is that the cone from a point over a sufficiently small neighborhood of the point is a linear space.

For a system of polynomials $F = (f_1,\dots,f_r)$, we have a containment 
$$
T(\VV(F)) \subset T(\VV(f_1)) \cap \dots \cap T(\VV(f_r)),
$$
which is strict, in general. The intersection on the right is called the \emph{tropical prevariety} of this system. 
It follows that a tropical variety $T(\VV(F))$ is $0$-dimensional, if the corresponding tropical prevariety can't contain a tropical curve.

For a polyhedral complex we define its \emph{recession cone} as the convex hull of the direction vectors for all halflines contained in the support of the complex together with the origin. It is straightforward to show the following.
\begin{proposition}\label{prop:pointed-recession-cone}
The recession cone of a polyhedral complex containing a tropical curve cannot be pointed.
\end{proposition}

\subsection{A fiber of tropicalization is Zariski dense}\label{subsec:dense}
The final ingredient that we need is the fact that the preimage of a point in the image of $T: (K^*)^n \to \QQ^n$ is a Zariski dense subset of $(K^*)^n$. A more general version of this fact is the result of~\cite{payne2009fibers}.
Note a formally redundant change from $N$ to $n$. We emphasize the need to use this fact on the space of masses.

\section{Computational method}\label{sec:computations}

We consider a collection of equations from~\Cref{sec:AC} as a system of polynomials $F = (f_1,\dots,f_r) \subset K[r_{12},\dots,r_{(n-1)n}]$ and specialize the masses $\masses\in K^*$.

By Kapranov's Theorem, one can obtain the tropical hypersurface $T(\VV(f_i))$ from a lifted Newton polytope of $f_i$ (\cite[Proposition~3.1.6]{tropicalbook}). Therefore, the task of computing the intersection $T(\VV(f_1)) \,\cap\, \dots \,\cap\, T(\VV(f_r))$ amounts to designing an efficient algorithm in polyhedral geometry. 

\begin{figure}[ht!]
\floatbox[{\capbeside\thisfloatsetup{capbesideposition={right,center},capbesidewidth=6.5cm}}]{figure}[\FBwidth]
{\caption{
A comet and its recession cone, which is pointed (i.e., contained in an open half-space). 
\label{fig:comet}}}
{\includegraphics[width=6.5cm]{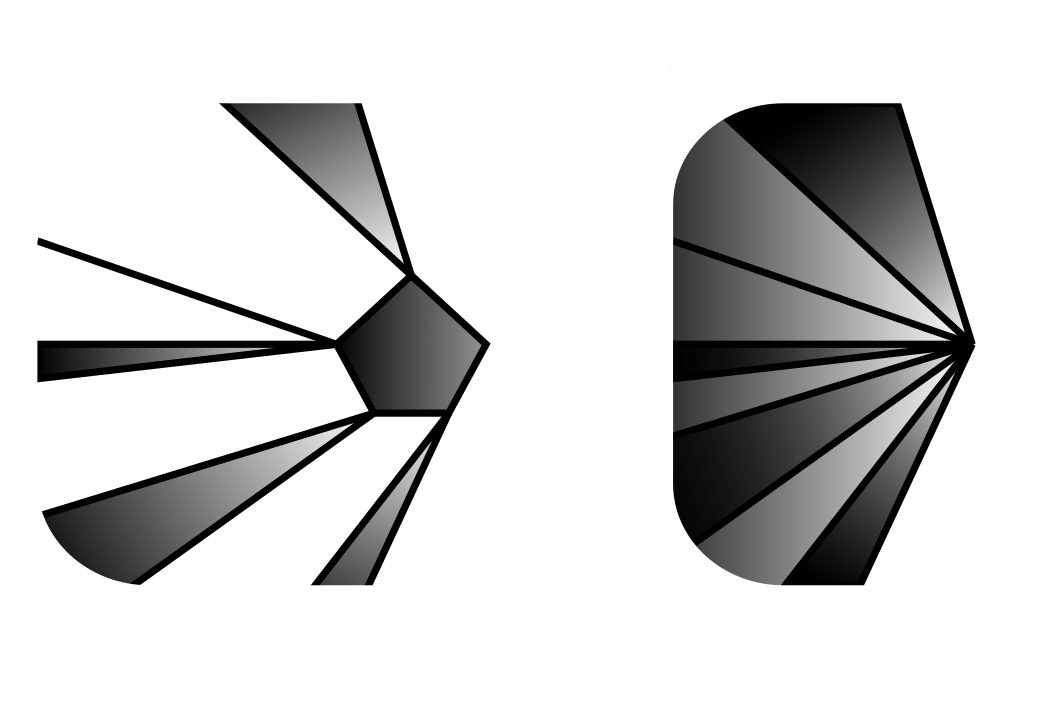}}
\end{figure}

Once this tropical prevariety is obtained, we can inspect its connected components --- we call these components \emph{comets} (see \Cref{fig:comet} for illustration\footnote{This is planar picture is only an illustration. The actual polyhedral complexes that we compute are in higher ambient dimension.}) --- for a possible containment of a tropical curve. 
{\bf If} this inspection yields a negative result in view of~\Cref{prop:pointed-recession-cone}, then we conclude that $T(\VV(F))$ is $0$-dimensional.

The polynomials $F$ that we use have coefficients that are linear in $\masses$. The following would be true more generally (for polynomial or rational-function expressions in masses): if we fix the valuations of masses and then consider only Puiseux series masses (with these valuations) such that there are no cancellations of leading terms, then the prevariety would depend only on these valuations (i.e., the image of the point $(\masses)$ under tropicalization  $T: (K^*)^n \to \QQ^n$). By Section~\ref{subsec:dense} the set of \emph{all} choices with fixed valuation is Zariski dense. 
This {\bf would} conclude the proof of~\Cref{conj:generic} since having $\VV(F)$ $0$-dimensional implies that $\VV(F)$ is finite.


\subsection{Proof for five bodies}
We consider the system of polynomials $F$ with
\begin{itemize}
    \item 10 symmetric Albouy-Chenciner equations
    \item 20 Albouy-Chenciner equations
    \item 5  Cayley-Menger determinants
\end{itemize}
The valuations of the masses are chosen to be $1,4,9,16$ and $25$, respectively. It is easiest to carry out the computation with software that is capable of dealing with polyhedral complexes that are not just fans. Version 0.7 of gfan~\cite{gfan} has this capability when computing tropical prevarieties.
The dynamic decomposition algorithm used and its implementation are described in~\cite{dynamicdecomposition}.
The software computes a polyhedral complex with the following f-vector:
$$(3586,12012,18531,15625,7072,1357)$$
For instance, the first number, $3586$, is the number of vertices in the complex. 

The complex has 257 connected components. 
For each component we collect its rays, which generate the recession cone. It turns out that all 257 are comets --- their recession cones are pointed --- making it impossible for the tropical prevariety to contain a balanced tropical curve. Hence, the variety defined by $F$ is zero-dimensional.

The time for computing the prevariety and checking pointedness is 80 cpu minutes (5 minutes on 16 threads) with 64-bit machine integers. Overflow checks are not always performed by the code. To justify that, one needs to consider that for our simplex algorithm implementation tableau entries are always subdeterminants of the starting tableau, providing a bound on the entries. Alternatively, we exploit the templated implementation and run the code with overflow checks in the integer class. This was done and takes about 10 times longer than the original computation. 

Different choices of valuations may either not work or make the proof easier. The argument fails for valuations $(0,1,2,3,4)$ because not all components have pointed recession. On the other hand powers of 3: $(1,3,9,27,81)$ give a prevariety with f-vector $(1506,4744,8586,8787,4652,993)$ with pointed recession cone, meaning that it is unnecessary to determine connected components.

\subsubsection{Verification in gfan}
We explain how the last claim can be verified in gfan version 0.7. The relevant polynomials are produced with
\begin{verbatim}
gfan _nbody -N5 --masses --alsosymmetric --cayleymenger2
\end{verbatim}
Because coefficients do not cancel it suffices to do the substitutions $m_i\mapsto t^{3^{i-1}}$, compute the tropical prevariety and study its unbounded directions. The substitution and change of ring is performed in the Linux shell with the one-line command
\begin{verbatim}
gfan _nbody -N5 --masses --alsosymmetric --cayleymenger2 |
sed -e"s/Q\[m1,m2,m3,m4,m5,/Q(t)\[/g" -e"s/m1/t1/g" 
-e"s/m2/t3/g" -e"s/m3/t9/g" -e"s/m4/t27/g" -e"s/m5/t81/g"
\end{verbatim}
The output is then used as input for
\begin{verbatim}
gfan _tropicalprevariety --usevaluation  -j16 --mint --minx
\end{verbatim}
which computes the tropical prevariety in 16 threads using the min-convention as in~\cite{tropicalbook}. To speed up computations and get reported timings \texttt{--bits64} can be used.

The output describes the fan over the complex $T(\VV(f_1)) \,\cap\, \dots \,\cap\, T(\VV(f_r))$. Notice that the support of this fan is not closed. Its lineality space is 0-dimensional and therefore the rays listed  should be interpreted either as points (those rays having first coordinates non-zero) or unbounded directions (those rays having first coordinates equal to zero).
By inspection it is now easy to see that among the unbounded directions all coordinates are non-positive, proving that the prevariety  has pointed recession cone.

With more tricks the computation can be verified in older versions of gfan too. Verifying other of our claims is more cumbersome and requires some programming.

\subsection{Special problems are harder} \label{sec:special}
For example, in case $n=4$, if all masses were chosen to have the same valuation, e.g., $(0,0,0,0)$, our tropical computation would be equivalent to that of~\cite{hampton2006finiteness} and, therefore, require further argument to reach the finiteness conclusion. The purely polyhedral computation using valuation $(0,0,0,0,0)$ leads to no conclusion for $n=5$.
Moreover, our experiments for specialization of the problem with a less restrictive symmetry, e.g., $m_1=m_2$ and $m_3=m_4$, setting valuations in the form $(a,a,b,b,c)$ fail to produce a proof either. Counterintuitively, there is more     hope for our methodology in the general case $n=6$ than in special cases of $n=5$.

\subsection{The six body case}
\label{sec:n=6}
For six bodies we were able to carry out the computation of the tropical prevariety for a system similar to the one for the five body case within 100 cpu days (about three days on 35 threads). We did this for a particular choice of mass valuations. Unfortunately the recession cones of the connected components were not all pointed. Therefore more experimentation with different valuations or different versions of the Albouy-Chenciner equations is needed. 

We are aware of another work, \cite{chang2023toward}, that attempts to conquer the $n=6$ case by developing symbolic computation algorithms following Albouy-Kaloshin's approach for $n=5$. For now, this attempt fails as well.

\section*{Acknowledgments} Although the main computational proof was carried out solely in gfan~\cite{gfan}, a software specializing in tropical geometry, a lot of preparatory work and experimentation was done in Macaulay2~\cite{M2}, a more general computer algebra system.
The research of the second author is partially supported by NSF DMS award 2001267 and the Simons Fellows program.

\bibliographystyle{abbrv}
\bibliography{refs}

\end{document}